\documentclass[reqno,11pt,centertags,draft]{amsart}
\date{\today}
\title[]{Remark on the paper "Asymptotic behavior of polynomials orthonormal on a homogeneous
set"}

\author[]{F. Peherstorfer,   and P. Yuditskii}

\begin{document}
\date{\today}

\maketitle

\begin{abstract}
Minor modifications are given to prove the Main Theorem under the
Blaschke (instead of Carleson) condition as well as  a small
historical comment. 
\end{abstract}

Because of the reference \cite{DKS} on our paper \cite{PY03} a
certain historical comment is needed.

In \cite{PY03} we generalized H. Widom's Theorem \cite{W} based on
an absolutely new idea, dealing with one dimensional perturbation of
a given Jacobi matrix. As it is well known if  a Jacobi matrix $J$
has the compact $E=[b_0,a_0]\setminus \cup_{j\ge 1}(a_j,b_j)$ as
spectral set, then its one dimensional perturbation may have in
addition to $E$ spectral points in the gaps (one point in each gap).
Since in our case, the set $E$ is possibly a Cantor type set, an
{\it infinite} number of spectral points $X$ has to be added to the
support of the spectral measure. A homogeneous set $E$ possesses the
following very nice property \cite{JM85}. Let
$B(z,z_0)=B(z,z_0;\Omega)$ be the Blaschke factor in the domain
$\Omega=\overline{ \mathbb C}\setminus E$ with zero at $z_0 \in
\Omega$. From each interval $(a_j,b_j)$ let us pick arbitrarily
exactly one $x_j$. Then
\begin{equation}\label{hr1}
    \inf_k\prod_{j\not =k}|B(x_k,x_j)|>0.
\end{equation}
Note that the convergence of the product $\prod_{j }|B(z_0,x_j)|>0$
(the Blaschke condition) corresponds to the so called Widom property
of the domain $\Omega$. So,  the domain with a homogeneous boundary
has even a better property: the more restrictive condition
\eqref{hr1}, the so called Carleson condition, holds for the given
Blaschke product.

Thus to use our idea on a one dimensional perturbation, we were
enforced to work with  spectral measures supported on a homogeneous
set $E$ but also having, possibly an infinite set of mass points
supported on an arbitrary (real) set $X$, satisfying (similar to
\eqref{hr1}) the Carleson condition
\begin{equation}\label{hr2}
    \inf_{x\in X}\prod_{y\in X, y\not =x}|B(x,y)|>0.
\end{equation}

When the results of \cite{PY03} were presented and the manuscript
was submitted for publication, we recognized that the math community
is much more interested in an infinite number of mass point than in
a Cantor type spectral sets. We wrote \cite{PY01} considering $E$
just as a single interval. It became, almost immediately clear that
the Carleson condition was too restrictive and a simple trick
\cite[(2.9)]{PY01} allows to use only Blaschke condition.

Unfortunately, the first paper  was finally published later than the
second one (in this way we have, formally, a more fresh paper with
the more restrictive condition \cite[(6.1)]{PY03}  on the set $X$).

{\it To prove  the Main Theorem under the Blaschke  condition} one
has to estimate \cite[(6.9)]{PY03} in the way \cite[(2.9)]{PY01}, of
course, following to the standard  strategy in this paper (Lemma
4.2, Lemma 5.3, etc):

1. For a fixed $\epsilon>0$ use the finite covering of $\Gamma^*$
(see \cite[p.139]{PY03}),
\begin{equation}\label{hr3}
\Gamma^*=\cup_{j=1}^{l(\epsilon)}\left\{\beta:
\text{dist}(\beta,\beta_j)\le\eta(\epsilon)\right\}.
\end{equation}

2. For the exhaustion $\{X_N\}$ of $X$ by finite sets
\cite[p.138]{PY03}, let $M_N$ be the character of the Blaschke
product $B_N$ with zeros at $X\setminus X_N$. Since $B_N(0)\to 1$
and $M_N\to 1_{\Gamma^*}$ we can choose $N$ so big that
\begin{equation}\label{hr4}
    1-B_N(0) \le \epsilon,   \ \  1-\Delta^{M_N^{-1}}(0) \le \epsilon
\end{equation}
(for the definition of $\Delta^{\alpha}$ see \cite[p.125]{PY03}).

3. Having a finite number of the reproducing kernels (characters
$\beta_j$) and a finite number of points ($X_N=\{\text{zeros of}\
B/B_N\}$), choose $n$ due to the  estimation
\begin{equation}\label{hr5}
\begin{split}
    &\left|\sum_{X_N}\left\{b^{n+1}K^{\beta_j}\left(\frac{-z'}{\phi'}\right)P_n\right\}\right|\\
    \le(\sup_{X_N}|b|)^n
    &\sqrt{\sum_{X_N}|P_n|^2\sigma_l}
    \sqrt{\sum_{X_N}\left
    |K^{\beta_j}\frac{b}{\psi}\right|^2\sigma_l}.
    \end{split}
\end{equation}

Let us mention that the key Lemmas 1.1, 2.2, 2.4, 5.2 were proved
under the Blaschke condition \cite[(2.5)]{PY03}.

\bibliographystyle{amsplain}

\end{document}